\documentclass[10pt, final]{amsproc}
\usepackage{amssymb}
\usepackage[mathscr]{eucal}
\usepackage{amscd}
\usepackage{amsmath}
\usepackage{setspace}
\usepackage{tikz}
\def\Ker{\operatorname{Ker}}

\def\dim{\operatorname{dim}}
\def\supp{\operatorname{supp}}

\theoremstyle{plain}
\newtheorem{theorem}{Theorem}
\newtheorem{problem}{Problem}
\newtheorem{proposition}{Proposition}
\newtheorem{lemma}{Lemma}
\newtheorem{corollary}{Corollary}
\newtheorem{definition}{Definition}
\theoremstyle{remark}

\title{A space with no unconditional basis that satisfies the Johnson-Lindenstrauss lemma}

\address{Instituto de Matem\'aticas Imuex, Universidad de Extremadura, Avenida de Elvas s/n, 06011 Badajoz, Spain.}
\email{jesus@unex.es}
\author{Jes\'us Su\'arez de la Fuente}

\subjclass[2010]{(Primary) 46B20, 46B06; (Secondary) 46B70, 46M18, 46B45}

\keywords{Johnson-Lindenstrauss lemma, Weak Hilbert, interpolation, twisted sum, centralizer}

\thanks{The author was supported in part by project MTM2016-76958-C2-1-P and project IB16056.}

\begin{document}

\begin{abstract}
We give the first example of a nontrivial twisted Hilbert space that satisfies the Johnson-Lindenstrauss lemma. This space has no unconditional basis. We also show that such a space gives a partial answer to a question of Mascioni.
\end{abstract}

\maketitle

\section{Introduction}


Let us begin by introducing the following definition as given in \cite{JoNa}:

\begin{definition} A Banach space $X$ satisfies the J-L lemma if given $x_1,...,x_n \in X$, there exists a linear mapping $L:X\to F$, where $F\subseteq X$ is a linear subspace of dimension $O(\log(n))$, such that for all $i, j \in \{1,...,n\}$ we have $$\|x_i-x_j\|\leq \|L(x_i)-L(x_j)\|\leq O(1)\|x_i-x_j\|.$$
\end{definition}

The definition is motivated by a celebrated result of Johnson and Lindenstrauss. They proved in \cite{JL} that Hilbert spaces do satisfy the J-L lemma. Their result has found many applications in mathematics and computer science. Due to the usefulness of the result there have been considerable effort by researches to prove it in non Hilbertian settings. This is, to find an $X$ different from a Hilbert space that still satisfies the sort of compression described above. The only positive result is due to Johnson and Naor \cite{JoNa} who showed that $T^2$, the $2$-convexification of the Tsirelson space, satisfies the J-L lemma. The space $T^2$ is the leading member of the class of weak Hilbert spaces introduced by Pisier in \cite{Pi}. And the result of Johnson-Naor can be extended to weak Hilbert spaces with an unconditional basis. However, it is not known if arbitrary weak Hilbert spaces satisfy the J-L lemma, see \cite[Remark 1]{JoNa}. A new weak Hilbert space called $Z(T^2)$ has been recently constructed in \cite{JS}. It appears as the derived space of $T^2$ with its dual, so it is a twisted Hilbert space. In particular, it must fail to have an unconditional basis by a result of Kalton. Despite this fact, we shall show that $Z(T^2)$ satisfies the J-L lemma. 

On the other hand, $Z(T^2)$ is related to a problem of V. Mascioni. Indeed, in \cite{M2} it was introduced a variant of the Homogeneous Banach space problem that reads as follows:
If $X$ is an infinite dimensional Banach space, and every infinite dimensional subspace of $X$ is isomorphic to its dual space, is $X$ isomorphic to a
Hilbert space? Motivated by this question, Mascioni studies the corresponding finite dimensional problem.
\begin{definition} A space is locally self dual (LSD) if there is a constant
$c>0$ such that every finite dimensional subspace of $X$ is $c$–isomorphic to its dual space.
\end{definition}

\begin{problem}\label{dos} Are LSD spaces isomorphic to Hilbert spaces?
\end{problem}

Mascioni then proves that LSD spaces are weak Hilbert spaces. So an hypothetical counterexample must be a weak Hilbert space.  Since twisted Hilbert spaces obtained by derivation are isomorphic to their own dual, it seems that $Z(T^2)$ could shed some light on Mascioni's questions. Indeed, we shall show that the space $Z(T^2)$ enjoys the property of being ``asymptotically" LSD in the following sense:
\begin{definition} A Banach space $X$ is asymptotically LSD if for every $\beta>1$  there is $c(\beta)>0$ such that given a finite dimensional subspace $E$ of $X$, we may find a further subspace $F$ of $X$ with $E\subseteq F$ for which $\dim F\leq \beta \dim E$ and such that $F$ is $c(\beta)$-isomorphic to $F^*$.
 \end{definition}

In this paper we use standard notation for Banach space theory as provided in \cite{AK}. To avoid unnecessary definitions and repetitions we shall assume the reader familiar with some background on complex interpolation and centralizers. The interested reader may find all necessary material on interpolation in the classic book of Bergh and L\"ofstr\"om \cite{BL} while for centralizers we follow the presentation given by Kalton and Montgomery-Smith in \cite{KM}. The interested reader may fulfill all the technical details in \cite{CCS}. We will only require a minimum knowledge on the class of weak Hilbert spaces; two good references being \cite{CS} and \cite{Pi}. Finally, a reference for the space $Z(T^2)$ is \cite{JS}.

\section{The Johnson-Lindenstrauss lemma}
Let us introduce some notation. `Space' means `infinite-dimensional Banach space' unless specified otherwise and $\omega$ denotes the vector space of all complex sequences. So let $X$ be a space with an unconditional basis $(e_j)_{j=1}^{\infty}$. Let us denote by $\mathcal F(X,X^*)$ the Calder\'on space \cite[Page 1138]{KM} and by $\delta_{1/2}:\mathcal F(X,X^*)\to (X,X^*)_{1/2}$ the natural evaluation map onto the interpolation space $(X,X^*)_{1/2}$ \cite[Page 1139]{KM}. Then we define the \textit{derived space}  $Z(X)$ of $(X,X^*)_{1/2}$ as 
$$Z(X)=\{ (x,y)\in \omega \times \omega: \|(x,y)\|_{Z(X)}<\infty \},$$ 

\noindent where $$\|(x,y)\|_{Z(X)}=\inf \{ \|F\|_{\mathcal F(X,X^*)}: F(1/2)=x, F'(1/2)=y\}.$$

In \cite{KM} the derived space is denoted $d(X,X^*)_{1/2}$. It is easy to check that $Z(X)$ is a twisted sum of $(X,X^*)_{1/2}$ \cite[Page 1159]{KM}; this is, it contains a subspace isomorphic to $(X,X^*)_{1/2}$ such that the quotient is $(X,X^*)_{1/2}$. Since it is well known that $(T^2,(T^2)^*)_{1/2}=\ell_2$ we find that $Z(T^2)$ is a twisted Hilbert space; we informally say that $Z(T^2)$ is obtained by derivation. The Kalton-Peck space $Z_2$ is also obtained by derivation and corresponds to $Z(\ell_p)$ for any $p\neq 2$. Let us give a word on centralizers. The natural norm on $Z(T^2)$ is hard to compute. However, we may use an equivalent quasi-norm. This quasi-norm is given in the following way:
$$\|(x,y)\|=\|x-\Omega(y)\|_2+\|y\|_2,$$
for some map $\Omega:\ell_2\longrightarrow \omega$, see \cite[Proposition 7.2.]{CCS}. A simple way to generate $\Omega$ is the following. Fix a constant $\rho>1$ and pick $B$ an homogeneous $\rho$-bounded selection for $\delta_{1/2}$. Then, it is not hard to check that $\Omega:=\delta_{1/2}'\circ B$ is such map, what in particular means that for some $C>0$
\begin{equation}\label{uncond}
\|\Omega(ax)-a\Omega(x)\|\leq C\|a\|_{\infty}\|x\|, \;\;a\in \ell_{\infty}, x\in \ell_2.
\end{equation}
We will write $\Omega_{\rho}$ instead of $\Omega$ if we need to emphasize the role of $\rho$. A map that satisfies the inequality above is called a \textit{centralizer}. More specifically, such $\Omega$ is a centralizer corresponding to $(T^2,(T^2)^*)_{1/2}=\ell_2$. There is of course some ambiguity in the choice of $\Omega$ but the difference of two choices is bounded, so they give equivalent quasi-norms. In general, the space $Z(X)$ has a natural basis that we denote as usual by $v_{2j-1}=(e_j,0)$ and $v_{2j}=(0,e_j)$ for all $j\in \mathbb N$.  In this way, $Z(X_n)$ or $[v_j]_{j=1}^{2n}$ will denote the closed linear span of $(v_j)_{j=1}^{2n}$. 

Recall that a deep result of Kalton \cite{ka2} shows that the only twisted Hilbert space (given by derivation or not) with an unconditional basis is $\ell_2$. So, in particular, $Z(T^2)$ must fail to have an unconditional basis. The key of the argument to show that $Z(T^2)$ satisfies the J-L lemma lies in the notion of $\log$-Hilbertian that Johnson and Naor introduced in \cite{JoNa} as follows:  We say that a space is $\log$-\textit{Hilbertian} provided that for every finite dimensional subspace $E$ of $X$ there are subspaces $E_1$ and $E_2$ of $X$ such that $E = E_1\oplus E_2$, $\dim(E_1) = O(\log(\dim E))$ and $d_{E_2} = O(1)$. 
 Johnson and Naor proved that if a space $X$ is $\log$-Hilbertian then $X$ satisfies the J-L lemma \cite[Section 3]{JoNa}. 

\begin{theorem}\label{log}
$Z(T^2)$ satisfies the J-L lemma.
\end{theorem}
\begin{proof} 
We only need to prove that $Z(T^2)$ is $\log$-Hilbertian. So given an $n$-dimensional subspace $E$ of $Z(T^2)$, write $E=E_1\oplus E_2$, where $E_1=E \cap [v_j]_{j=1}^{[\log(n)]}$ and $E_2=E \cap [v_j]_{j=[\log(n)]+1}^{\infty}$ (we have denoted $[\log(n)]$ the integer part). By \cite[Theorem 1, Step I]{JS}, there is $C>0$ such that every $n$-dimensional subspace of $[v_j]_{j=[\log(n)]}^{\infty}$ is $C$-isomorphic to Hilbert, so in particular $d_{E_2}\leq C$.
\end{proof}

\section{On Mascioni's question}

 One of the key arguments to prove our partial answer to Mascioni's question is given by the following result on twisted sums. Since it is hard for a non-specialist to extract it from the existing literature, we include a detailed proof for the sake of completeness. Let $X$ denote a space with an unconditional basis $(e_j)_{j=1}^{\infty}$ and let $Z(X)$ be the twisted Hilbert space given by derivation as in the previous section, then

\begin{lemma}\label{isomorphism} There is $\lambda>0$ such that for all $n=1,2,...$ the space $Z(X_n)$ is $\lambda$-isomorphic to its dual. 
\end{lemma}
\begin{proof}
Fix $\rho>1$ and let $\Omega_{\rho}$ be a centralizer corresponding to $(X,X^*)_{1/2}$. This means, in particular, that $\Omega_{\rho}=\delta_{1/2}'\circ B$, where $B$ is a homogeneous $\rho$-bounded selector for $\delta_{1/2}$. Let us denote by $[e_j]_{j=1}^{n}$ and $[e_j^*]_{j=1}^{n}$ the corresponding linear span of the elements in the canonical basis of $X$ and $X^*$ respectively, and by $\Omega_{\rho}^n$ the restriction of $\Omega_{\rho}$ to $([e_j]_{j=1}^{n},[e_j^*]_{j=1}^{n})_{1/2}$. Since $[v_j]_{j=1}^{2n}$ and $d_{\Omega_{\rho}^n}([e_j]_{j=1}^{n},[e_j^*]_{j=1}^{n})_{1/2}$ have, up to an absolute constant, $\rho$-equivalent norms by \cite[Proposition 7.2.]{CCS}, we shall prove that the isomorphism constant of  $d_{\Omega_{\rho}}([e_j]_{j=1}^{n},[e_j^*]_{j=1}^{n})_{1/2}$ with its dual depends only of $\Delta(\Omega_{\rho})$ for all $n=1,2,...$, which will be clearly enough. Recall that by $\Delta(\Omega_{\rho})$ we denote as usual the least constant $C>0$ such that
\begin{equation*}
\|\Omega_{\rho}(ax)-a\Omega_{\rho}(x)\|\leq C\|a\|_{\infty}\|x\|,
\end{equation*}
that is guaranteed by the very definition of centralizer. 
 So let us begin observing that $d_{\Omega_{\rho}^n}([e_j]_{j=1}^{n},[e_j^*]_{j=1}^{n}))_{1/2}$ is isometric to $d_{-\Omega_{\rho}^n}([e_j]_{j=1}^{n},[e_j^*]_{j=1}^{n}))_{1/2}$ via the map $(x,y)\to (-x,y)$. We shall show that the duality map $$I:d_{\Omega_{\rho}^n}([e_j]_{j=1}^{n},[e_j^*]_{j=1}^{n})_{1/2}\longrightarrow \left(d_{-\Omega_{\rho}^n}([e_j]_{j=1}^{n},[e_j^*]_{j=1}^{n})_{1/2}\right)^*$$ defined by $\langle I(x,y), (u,v) \rangle=\langle x,v \rangle+ \langle y,u \rangle$ is an onto isomorphism with isomorphism constant depending only of $\Delta(\Omega_{\rho})$. A quick check will show that we have a commutative diagram
$$
\begin{CD}
0 @>>>\ell_2@>j>>d_{\Omega_{\rho}^n}\ell_2@>q>>\ell_2@>>>0\;\;\;\;\\
&& @ViVV  @VIVV @ViVV\\
0 @>>>\ell_2^*@>q*>>\left(d_{-\Omega_{\rho}^n}\ell_2\right)^*@>j*>>\ell_2^*@>>>0,\;\;\;\;\\
\end{CD}
$$
where $j(x)=(x,0), q(x,y)=y$ and $i(x)(v)=\langle x,v \rangle$. Thus, it follows easily that $I$ is injective and surjective, see the 3-lemma in \cite{CG}, since $i$ is an isomorphism. Let us prove that $I$ is bounded:
\begin{eqnarray*}
|\langle x,v \rangle+ \langle y,u \rangle|&=& |\langle x-\Omega^n_{\rho}(y),v \rangle+ \langle y,u+\Omega^n_{\rho}(v) \rangle +\langle \Omega^n_{\rho}(y),v\rangle + \langle y,-\Omega^n_{\rho}(v)\rangle|\\
&\leq& \|x-\Omega^n_{\rho}(y)\|\|v\|+\|y\|\|u+\Omega^n_{\rho}(v)\|+8\Delta(\Omega_{\rho})\|y\|\|v\|\\
&\leq& 8\Delta(\Omega_{\rho})\|(x,y)\|\|(u,v)\|,
\end{eqnarray*}
where the first inequality follows by $\Delta(\Omega_{\rho}^n)\leq\Delta(\Omega_{\rho})$ and the bound $$|\langle \Omega^n_{\rho}(y),v\rangle - \langle y,\Omega^n_{\rho}(v)\rangle|\leq 8\Delta(\Omega_{\rho}^n)\|y\|\|v\|,$$
to be found in the discussion before \cite[Corollary 4]{FC}; it also follows from (5.5) in \cite[Theorem 5.1.]{kaltmem} taking $f_1=g_2=y$ and $f_2=g_1=v$. 

We use now an argument of \cite[Proposition 16.12.]{BeLi} to keep track of the norm of the inverse. So fix $(x,y)$ and pick $v$ with $\|v\|=1$ such that $\langle x-\Omega^n_{\rho}(y),v \rangle=\|x-\Omega^n_{\rho}(y)\|$ and let $u=-\Omega_{\rho}^n(v)+ (8\Delta(\Omega_{\rho})+1)y/\|y\|$. It is trivial that $u+\Omega_{\rho}^n(v)\in \ell_2$ and thus $(u,v)\in d_{-\Omega_{\rho}^n}([e_j]_{j=1}^{n},[e_j^*]_{j=1}^{n})_{1/2}$. Therefore, we find
\begin{eqnarray*}
|\langle x,v \rangle+ \langle y,u \rangle|&=& |\|x-\Omega^n_{\rho}(y)\|+ (8\Delta(\Omega_{\rho})+1)\|y\|+ \langle \Omega^n_{\rho}(y),v\rangle + \langle y,-\Omega^n_{\rho}(v)\rangle|\\
&\geq&\|x-\Omega^n_{\rho}(y)\|+ (8\Delta(\Omega_{\rho})+1)\|y\|-8\Delta(\Omega_{\rho})\|y\|\\
&\geq& \|(x,y)\|.
\end{eqnarray*}
Since $\|(u,v)\|=8\Delta(\Omega_{\rho})+2$, we finally find that the isomorphism constant of $d_{\Omega_{\rho}^n}([e_j]_{j=1}^{n},[e_j^*]_{j=1}^{n})_{1/2}$ with its dual depends only of $\Delta(\Omega_{\rho})$. It is not hard to prove the estimate $\Delta(\Omega_{\rho})\leq 4\rho$ once the bound $\| \delta_{\theta}':\Ker \delta_{\theta}\to X_{\theta}\|\leq \max\{ \theta^{-1},(1-\theta)^{-1}\}$ is provided, see \cite[Page 4676]{CFG}. Since $d_{\Omega_{\rho}^n}([e_j]_{j=1}^{n},[e_j^*]_{j=1}^{n})_{1/2}$ and $[v_j]_{j=1}^{2n}$ have, up to an absolute constant, $\rho$-equivalent norms by \cite[Proposition 7.2.]{CCS}, we find that $[v_j]_{j=1}^{2n}$ is $\lambda$-isomorphic to its dual with $\lambda=f(\rho)$ for some fixed function $f$.
\end{proof}
We are ready now to prove our partial answer to Problem \ref{dos}.
\begin{theorem}\label{mascioni}
$Z(T^2)$ is asymptotically LSD.
\end{theorem}
\begin{proof} We shall prove a slightly stronger result, that is, for every $m=1,2,...$ there is $c_m>0$ such that for every $n$-dimensional subspace $E$ of $Z(T^2)$ there is a subspace $F_m$ of $Z(T^2)$ verifying:
\begin{enumerate}
\item $E\subseteq F_m$.
\item $\dim F_m\leq n+\log_m(n)$.
\item $F_m$ is $c_m$-isomorphic to $(F_m)^*$. 
\end{enumerate}
The claim of the theorem follows trivially from this. So fix $m\in \mathbb N$ and pick $n_0=n_0(m)\in \mathbb N$ such that for all $n\geq n_0$ we have $$1\leq [\log_m(n)]-1,$$ where $[\log_m(n)]$ denotes the integer part. Let $E$ be a $n$-dimensional subspace of $Z(T^2)=[v_j]_{j=1}^{\infty}$ and assume first that $n\leq n_0$. Since by John's theorem \cite{JT} we always have the trivial bound $d(E,E^*)\leq \dim E$, we have the theorem proved with say $F_m=E$ and $c_m=n_0$. So we may assume that $n\geq n_0$. If $[\log_m(n)]$ is even, then write $E=E_1\oplus E_2$, where $E_1=E \cap [v_j]_{j=1}^{[\log_m(n)]}$ and $E_2=E \cap [v_j]_{j=[\log_m(n)]+1}^{\infty}$. By \cite[Theorem 1, Step I]{JS}, there is $C>0$ such that every $n$-dimensional subspace of $[v_j]_{j=[\log_m(n)]}^{\infty}$ is $C^m$-isomorphic to Hilbert, so in particular $d_{E_2}\leq C^m$. By Lemma \ref{isomorphism},  $[v_j]_{j=1}^{2n}$ is $\lambda$-isomorphic to its dual for every $n\in \mathbb N$. Then, it is a straightforward computation to check that $F_m:=[v_j]_{j=1}^{[\log_m(n)]}\oplus E_2$ is  $f(C^m)$-isomorphic to its dual  for some fixed function $f$. Therefore, such $F_m$ satisfies the claim of the proposition. Indeed, we have just seen $(3)$ and $(1)$ is clear by construction. For $(2)$ just observe that $$\dim F_m=[\log_m(n)]+\dim E_2\leq [\log_m(n)]+n.$$  If $[\log_m(n)]$ is odd, then consider the splitting above with $E_1=E \cap [v_j]_{j=1}^{[\log_m(n)]-1}, E_2=E \cap [v_j]_{j=[\log_m(n)]}^{\infty}$ and argue as before. Finally, take $$c_m:=\max \{ f(C^m),  n_0(m) \}.$$
\end{proof}

\section{Further considerations on $Z(T^2)$}\label{last}
Motivated by Theorem \ref{mascioni}, our goal in this section is to study more in detail the basis of $Z(T^2)$.  In this line, we complete \cite[Corollary 1]{JS}, where it was shown that $d_{Z(T^2_n)}=o(\log_m(n))$ for $m=1,2,...$, by checking that the equivalence constant of the basis of $Z(T^2)$ with the canonical one in $\ell_2$ is of the same order of magnitude as in $T^2$. We also introduce a parameter $D_n(X)$ as a substitute for the property $(H)$ on blocks basis sequences in a space $X$. We show that $D_n(Z(T^2))$ is of the order $\sqrt n$ and, as a consequence, we find an alternative and elementary proof that the Kalton-Peck space $Z_2$ cannot be a subspace of $Z(T^2)$, see the remark after \cite[Proposition 4]{JS}. Summing up, all the results in this second part show that, even if the basis of $Z(T^2)$ is not unconditional, it looks very much like the basis of $T^2$. We first deal with the finite-dimensional sections of the basis while separately in Subsection \ref{bloques} we give our estimate for block basis sequences.

Thus, let us give a few words on $Z(T^2_n)$ that is a twisted sum of $\ell_2^n$.
These twisted sums are not uniformly trivial, meaning that the distance to the corresponding Hilbert copy grows to infinity as the dimension tends to infinity.
\begin{proposition}\label{going}
$$\sup_{n\in \mathbb N} d_{Z(T^2_n)}=\infty$$
\end{proposition}
\begin{proof}
Assume on the contrary that $\sup_{n\in \mathbb N} d_{Z(T^2_n)}=K<\infty,$ so that $Z(T^2_n)$ has Rademacher type $2$ constant uniformly bounded by $K$.  If $\Omega$ denotes a centralizer corresponding to $Z(T^2)$, then we may assume that the restriction of $\Omega$ to $[e_j]_{j=1}^n$ is a centralizer corresponding to $Z(T^2_n)$. Therefore, we have for each $n\in\mathbb N$
\begin{equation}\label{estimate}
\mathbb E \left\|\left( 0,\sum_{j=1}^n r_j\lambda_je_j\right) \right\|_{Z(T^2)}\leq K \left(\sum_{j=1}^n |\lambda_j|^2 \right)^{1/2}.
\end{equation}

It follows readily from the definition of centralizer by picking $a\in\{-1,1\}$ in (\ref{uncond}) that the sequence $(v_{2j})_{j=1}^{\infty}$ is $\lambda$-unconditional for some $\lambda>0$. Thus, for every $n\in \mathbb N$

$$\left\|\Omega\left( \sum_{j=1}^n \lambda_je_j\right)\right\|+ \left\|\left( \sum_{j=1}^n \lambda_je_j\right)\right\|= \left \|\left( 0,\sum_{j=1}^n \lambda_je_j\right)\right\| \approx \mathbb E \left\|\left( 0,\sum_{j=1}^n r_j\lambda_je_j\right)\right\|.$$

In particular, combining this last equation with (\ref{estimate}), we find that $\Omega$ is bounded by a constant depending only on $K$. Another way to state this claim is the following: $\Omega$ is equivalent to the centralizer $0$; meaning that the difference takes values in $\ell_2$ and  is bounded. Observe that the centralizer $0$ is the one that corresponds to $(\ell_2,\ell_2)_{1/2}=\ell_2$. Thus, we are in position to apply Kalton's uniqueness theorem \cite[Theorem 7.6.]{ka} with $X_0=T^2$, $X_1=(T^2)^*$ and $Y_0=Y_1=\ell_2$. Then we have $T^2\approx \ell_2$, what is absurd.
\end{proof}
We prove now the announced result:
\begin{proposition}\label{fd}
There is $K>0$ such that for all $N=1,2,...$, $(v_j)_{j=1}^{N}$ is $K^{m}\log_m(N)$-equivalent to the unit vector basis of $\ell_2^{N}$ where $m=1,2,...$  
\end{proposition}
\begin{proof}
Let us assume first that $N=2n$ and recall that the constant of equivalence of $(e_j)_{j=1}^n$ in $T^2$ with $\ell_2^n$ is of the order $K^m\log_m(n)$, for some $K>0$, what follows from the proof of \cite[Theorem IV.b.3.]{CS}, see also \cite[Proposition 2.8]{NC}. This means in particular that 
$$(K^m\log_m(n))^{-1}\|y\|_{\ell_2}\leq \|y\|_{T^2}  \leq \|y\|_{\ell_2},\;\;y\in [e_j]_{j=1}^n.$$
Therefore, by simple duality 
$$\max \left( \|y\|_{T^2}, \|y\|_{(T^2)^*}\right)\leq K^m\log_m(n)\|y\|_2,\;\;y\in [e_j]_{j=1}^n.$$
For a fixed $y\in [e_j]_{j=1}^n$ and every $z\in \{ \omega: 0\leq \Re w\leq 1 \}$ consider the constant function in the space of Calder\'on $$B(y)(z)=y\in \mathcal F(T^2,(T^2)^*).$$ Therefore $B(y)$ is a $K^m\log_m(n)$-bounded selection for $y$. Moreover, we find $$\Omega_{K^m\log_m(n)}(y)=0,\;\;\;y\in [e_j]_{j=1}^n.$$ So that given $(x,y)\in[v_j]_{j=1}^{2n}$ we have that
\begin{equation}\label{pum}
\left \| \left (x, y\right) \right \|_{d_{\Omega_{K^m\log_m(n)}\ell_2}} =\left \|x-\Omega_{K^m\log_m(n)}(y) \right \| + \|y \|=\left \|x\right \|+\left \|y \right \|.
\end{equation}
The norm in $Z(T^2)$ is, up to an absolute constant,  $K^m\log_m(n)$-equivalent to the quasi-norm (\ref{pum}) (see \cite[Proposition 7.2.]{CCS}). Since (\ref{pum}) is clearly $\sqrt{2}$-equivalent to the norm in $\ell_2^{2n}$, the claim follows if $N=2n$ and thus trivially also if $N$ is odd.
\end{proof}
\subsection{An estimate for block basis sequences in $Z(T^2)$}\label{bloques}
Let us begin by recalling that a space $X$ has the \textit{property $(H)$} if for every $\lambda>1$ there is a constant $f(\lambda)$ such that for any normalized $\lambda$-unconditional basic sequence $x_1,...,x_n$ we have 
\begin{equation}\label{Hache}
f(\lambda)^{-1}\sqrt{n}\leq \left \|\sum_{j=1}^n x_j\right \|\leq f(\lambda)\sqrt{n}.
\end{equation}

Since it is well known that weak Hilbert spaces have the property $(H)$ \cite[Proposition 14.2.]{Pi}, $Z(T^2)$ has the property $(H)$ by \cite[Theorem 1]{JS}. However, the basis $(v_j)_{j=1}^{\infty}$ is not unconditional and thus we have not an estimate like (\ref{Hache}) for successive normalized blocks of the basis in $Z(T^2)$. We prove in Proposition \ref{parameterD} a soft version of the property $(H)$ for such blocks. 

Given two blocks $u_1,u_2$ of a basis, we write $u_1<u_2$ as usual for $\max u_1<\min u_2$. We also write $u_1\ll u_2$ for $1+\max u_1<\min u_2$, for which we observe the following.
\begin{lemma}\label{separation}
Let $u_j=(a_j,b_j)$ be blocks in $Z(T^2)$ for $j\leq 2$. Then
  \[ u_1\ll u_2 \Longrightarrow\left\{ \begin{array}{ll}
      \displaystyle a_1-\Omega(b_1)< a_2-\Omega(b_2),\\
      b_1<b_2, 
        \end{array} \right. \] 
for any centralizer $\Omega$ such that $\supp  \;\Omega(x)\subseteq \supp (x)$ for every $x\in \ell_2$.
\end{lemma}
\begin{proof}
Write $u_1=\sum_{j=n_0}^{n_1} \lambda_jv_j$, let $n_1=2j_0-1$ and thus $v_{n_1}=(e_{j_0},0)$. 
Since $\min u_2\geq n_1+2=2j_0+1$ we find
$$u_1\in [(e_1,0),(0,e_1),...,(e_{j_0},0)],$$
$$u_2\in [(e_{j_0+1},0),(0,e_{j_0+1}),(e_{j_0+2},0),...].$$
And thus the claim follows trivially. 
The case $n_1=2j_0$ is simpler.
\end{proof}
In what follows, let $X$ be a space with basis and let us write $u_1,...,u_n$ for $n$ blocks of such basis. We introduce the following parameter that will be our substitute for the property $(H)$:
 $$D_n(X)=\sup  \left \{    \|u_1+...+u_n\|_X :u_1<...< u_n, \|u_j\|\leq 1\right\}.$$
First, for the sake of completeness, we prove here
\begin{lemma}\label{symm} There are $c_1,c_2>0$ such that
 $$c_1^{-1}\sqrt{n}\leq D_n(T^2)\leq c_1\sqrt{n},\;\:n\in \mathbb N,$$
 $$c_2^{-1}\sqrt{n}\leq D_n((T^2)^*)\leq c_2\sqrt{n},\;\:n\in \mathbb N.$$
\end{lemma}
\begin{proof}
Let us deal with the case of $D_n(T^2)$. The lower bound is an easy consequence of the fact that for any admissible set $A\subseteq \mathbb N$ we have that $2^{-1}\sqrt{|A|}\leq \left \| \sum_{j\in A} e_j \right\|_{T^2}$. For the upper bound observe that, by unconditionality, the supremum must be attained on $n$ successive blocks in the unit sphere. Since $T^2$ has the property $(H)$ we are done. A similar argument works for $(T^2)^*$.
\end{proof}
We prove now our soft version of the property $(H)$ for blocks. 
\begin{proposition}\label{parameterD}
There is $C>0$ such that
 $$C^{-1}\sqrt{n}\leq D_n(Z(T^2))\leq C\sqrt{n},\;\:n\in \mathbb N.$$
\end{proposition}
\begin{proof}
It is enough to prove the claim in the quasi-norm induced by a centralizer. The lower bound is trivial by taking $v_1=(e_1,0)\ll v_3=(e_2,0)\ll...\ll v_{2n-1}=(e_n,0)$. For the upper bound, there is no loss of generality assuming that the centralizer satisfies $\supp  \;\Omega(x)\subseteq \supp (x)$ for every $x\in \ell_2$; see, for example, the remark after \cite[Corollary 1]{CSC}. We do first the case in which the blocks are `well separated'. Thus, pick any $n$ blocks $u_j=(a_j,b_j)$ in $Z(T^2)$ such that $u_1\ll...\ll u_n$ and 
\begin{equation}\label{normauno}
\|a_j-\Omega(b_j)\|+\|b_j\|=\|u_j\|\leq 1,\;\;j\leq n.
\end{equation} 
Then 
\begin{eqnarray*}
\left \| \sum_{j=1}^n u_j \right \|  &\leq & \left \| \sum_{j=1}^n (a_j - \Omega(b_j))\right\| + \left \|  \sum_{j=1}^n b_j \right \|+ \left \| \Omega\left(\sum_{j=1}^n b_j \right)- \sum_{j=1}^n\Omega (b_j)\right\|.\\
\end{eqnarray*}
Since $u_j\ll u_{j+1}$ holds for $j\leq n-1$, we find by Lemma \ref{separation} that $a_j-\Omega(b_j)< a_{j+1}-\Omega(b_{j+1})$ and $b_{j}<b_{j+1}$ for $j\leq n-1$. Thus by (\ref{normauno})
$$\left \| \sum_{j=1}^n u_j \right \| \leq 2\sqrt{n}+  \left \| \Omega \left(\sum_{j=1}^n b_j \right)- \sum_{j=1}^n\Omega (b_j)\right\|.$$
We need to bound the remaining term and we do it using the inequality provided by Castillo, Ferenczi and Gonz\'alez in \cite[Lemma 4.8.]{CFG} for the choice $M_{X,\mathcal C}(n)=D_n(X)$. Observe that we may do it since $b_{j}<b_{j+1}$ for $j\leq n-1$. Therefore we find:
$$\left \| \Omega \left(\sum_{j=1}^n b_j \right)- \sum_{j=1}^n\Omega (b_j)\right\|\leq 6D_n(T^2)^{1/2}D_n((T^2)^*)^{1/2}+ \left \|  \log \frac{D_n(T^2)}{D_n((T^2)^*)} \sum_{j=1}^n b_j\right \|.$$
By Lemma \ref{symm} and using again that $b_{j}<b_{j+1}$ for $j\leq n-1$ we are clearly done. For the general case pick $u_1<...<u_n$ with $\|u_j\|\leq 1$ for $j\leq n$ and suppose first that $n=2k$. Since $u_{2j-1} \ll u_{2j+1}$ and $u_{2j} \ll u_{2j+2}$ for every index $j$, we trivially have
$$\left \|\sum_{j=1}^{2k} u_j \right \|_{Z(T^2)} \leq \left \|\sum_{j=1}^{k} u_{2j-1} \right \|_{Z(T^2)} + \left \|\sum_{j=1}^{k} u_{2j} \right \|_{Z(T^2)}\leq c\sqrt{k},$$
where the last inequality holds by our previous discussion. 
The case $n=2k-1$ is similar.
\end{proof}
Since $Z(T^2)$ is a weak Hilbert space and $Z_2$ is not weak Hilbert by a result of Mascioni \cite[Proposition 4.1.]{M}, it follows that the Kalton-Peck space $Z_2$ is not isomorphic even to a subspace of $Z(T^2)$. We may, however, give an elementary proof of this fact. Let us briefly recall the definition of the generalized Kalton-Peck spaces \cite{kaltpeck}. Given a Lipschitz map $\phi:[0,\infty)\to \mathbb R$ with $\phi(0)=0$, the space $Z_2(\phi)$ is the twisted Hilbert space $$\{(x,y)\in \omega \times \omega: \|x-\Omega(y)\|_2+\|y\|_2<\infty\},$$ where $\Omega$ is the Kalton-Peck centralizer:
\begin{eqnarray}\label{kaltpeck2}
\Omega(y)=\sum_{j=1}^{\infty} y_j \phi \left(\log \left( \|y\|_2/|y_j|\right)\right)e_j,
\end{eqnarray}
taking $0\cdot \phi(\log \infty)$ as $0$ and $y=\sum y_je_j$. The case $\phi(t)=t$ corresponds to the so-called space $Z_2$.
\begin{corollary}\label{alternative} $Z_2(\phi)$ is not isomorphic to a subspace of $Z(T^2)$ for any unbounded $\phi$.
\end{corollary}
\begin{proof}
Suppose $Z_2(\phi)$ embeds isomorphically into $Z(T^2)$. Therefore, by a standard gliding hump argument (\cite[Proposition 1.3.10]{AK}), there exists a subsequence  $(v_{2k_j})_{j=1}^{\infty}$ in $Z_2(\phi)$ that is $K$-equivalent to a block basic sequence $(u_j)_{j=1}^{\infty}$ in $Z(T^2)$ with $u_j<u_{j+1}$ for every $j\in \mathbb N$. Thus, for the vectors $v_{2k_j}=(0,e_{k_j})$ and each $n\in \mathbb N$ we have
\begin{eqnarray*}
\sqrt{n}\phi(\log \sqrt{n})+\sqrt{n}=\left \| \sum_{j=1}^{n} (0,e_{k_j}) \right\|_{Z_2(\phi)}&\leq&K \left \|  \sum_{j=1}^{n} u_{j} \right \|_{Z(T^2)}\leq  KC \sqrt{n},
\end{eqnarray*}
where the equality holds by the definition of the quasi-norm in $Z_2(\phi)$ and the last bound is due to Proposition \ref{parameterD} for some absolute constant $C>0$. This last cannot hold for $n$ large enough if $\phi$ is unbounded what finishes the proof.
\end{proof}
\subsubsection{Comparison with other $Z(T^p)$}
Let us denote by $T^p$ the $p$-convexification of the Tsirelson space and by $(T^p)^*$ its dual for $1\leq p <\infty$. Using \cite[Proposition 6.1.]{CFG} and the unconditionality of the basis in $(T^p)^*$, we find that $(T^p,(T^p)^*)_{\frac{1}{2}}=\ell_2$ with equal norms. Therefore, we may also consider the twisted Hilbert space $$Z(T^p)=d(T^p,(T^p)^*)_{\frac{1}{2}},$$
for every $1\leq p <\infty$ with $p\neq 2$. Let us prove that we may also distinguish all these twisted Hilbert spaces from $Z(T^2)$ using only Proposition \ref{parameterD}. 
\begin{proposition}\label{Tsirelp}
	If $Z(T^p)$ embeds isomorphically into $Z(T^2)$ then $p=2$.
\end{proposition}
\begin{proof}
Assume we have such an embedding and find a subsequence in $Z(T^p)$, say $\{(0,e_{j}): j\in \mathbb N_1\}$ with $\mathbb N_1\subseteq \mathbb N$, that is $K$-equivalent to a normalized block basis sequence $\{u_j: j\in \mathbb N_1\}$ in $Z(T^2)$. In particular, for each admissible set $A\subseteq \mathbb N_1$:
	\begin{equation}\label{cuenta}
	\left \| \Omega\left ( \sum_{j\in A}e_j  \right)  \right\|+ \left\| \sum_{j\in A}e_j   \right\|=\left \| \sum_{j\in A} e_{j} \right\|_{Z(T^p)}\leq K \left \|  \sum_{j\in A} u_j \right \|_{Z(T^2)}\leq K' \sqrt{|A|},
	\end{equation}
	
	where the last inequality holds by Proposition \ref{parameterD}. Let us check that $\Omega$ is a multiple of the Kalton-Peck map \cite{kaltpeck}. Indeed, pick any vector $x$ whose support, say $A$, is admissible and such that $\|x\|_{\ell_2}=1$. Observe that we have $$\frac{1}{2^{1/p}} \left( \sum_{j\in A} |x_j|^{2}\right)^{\frac{1}{p}}\leq \| |x|^{\frac{2}{p}}\|_{T^p} \leq \left( \sum_{j\in A} |x_j|^{2}\right)^{\frac{1}{p}}=1,$$ and by simple duality, where as usual $\frac{1}{p}+\frac{1}{q}=1$, we obtain $$\| |x|^{\frac{2}{q}}\|_{(T^p)^*} \leq 2^{1/p} \left( \sum_{j\in A} |x_j|^{2}\right)^{\frac{1}{q}}=2^{1/p}.$$
	
	Since $$|x|=\left(|x|^{\frac{2}{p}}\right)^{\frac{1}{2}}\left(|x|^{\frac{2}{q}}\right)^{\frac{1}{2}},$$
	
	the expression 
	$$\frac{x}{|x|}\left(|x|^{\frac{2}{p}}\right)^{1-z}\left(|x|^{\frac{2}{q}}\right)^{z},$$
	defines a $2^{1/p}$-bounded selection at $x$ for the quotient map $\delta_{\frac{1}{2}}$. Thus, derivating and evaluating at $1/2$, we find:
	
	$$\Omega(x)=x\log \left( \frac{  {|x|^{\frac{2}{q}}}} {|x|^{\frac{2}{p}}}   \right) =\left( \frac{2}{q}-\frac{2}{p}\right) x \log |x|,$$
	
	as claimed. If we extend by homogeinity the expression above and use it in (\ref{cuenta}), we have:
	
	$$\left| \frac{2}{q}-\frac{2}{p}\right|\log \sqrt{|A|}\leq K'.$$
	
	But this only holds for arbitrary large $|A|$ if $p=q=2$ and thus we are done. 
\end{proof}
The reader may be tempted to believe that $Z(T^p)$ may be still a weak Hilbert space also for $p\neq 2$. However, this is not in the cards. The reason is that we have implicitly proved above that $Z(T^p)$ contains uniformly isomorphic copies of $Z_2^n$ for every $n\in \mathbb N$ whenever $p\neq 2$. This fact and a result of Bourgain, Tzafriri and Kalton about weak cotype 2 which we show in Proposition \ref{wc2} is enough to reach our conclusion. Recall that if we give for each $0<\delta<1$, the quantity (\cite{Pi}) $$d_X(\delta)=\sup_{E\subseteq X, \dim E <\infty }\inf \{ d(F,\ell_2^{dim F}): F\subseteq E, \dim F\geq \delta \dim E  \},$$
 then we say $X$ has \textit{weak cotype 2} if (and only if) for each $0<\delta<1$, $d_X(\delta)$ is finite. And just as a reminder, $X$ is a weak Hilbert space if and only if $X$ has weak cotype 2 and weak type 2.
\begin{proposition}\label{wc2}
Let $X$ be a space containing uniformly $(Z_2^n)_{n=1}^{\infty}$. Then $X$ fails to have weak cotype 2. 
\end{proposition}
\begin{proof}
Let us assume that $X$ has weak cotype $2$. Thus, given $0<\delta<1$, $Z_2^n$ should contain a $d_X(\delta)$-isomorphic copy of $$\ell_2^{\delta\cdot 2n}.$$ However, by the result of Bourgain, Tzafriri and Kalton \cite[Theorem 1.5]{BTK}, the only isomorphic copies (with isomorphism constant independent of the dimension) of a Hilbert space inside $Z_2^n$ have dimension bounded by $n +o(n)$. But then, for $\frac{1}{2}<\delta<1$, the constant $d_X(\delta)$ must depends on $n$ since otherwise $$\delta\cdot 2n \leq n + o(n),$$ and thus  $$(2\delta - 1 )n \leq o(n)$$ with $(2\delta - 1)>0$. And this last is absurd because $o(n)/n\longrightarrow 0$.
\end{proof}


Incidentally, we may prove that there is $C>0$ such that $Z_2^{n}$ contains $C$-isomorphic (and explicit) copies of $\ell_2^{n+ [\log n]}$ for every $n\in \mathbb N$, as theoretically allowed by the result of Bourgain, Tzafriri and Kalton. The proof of this fact follows easily localizing the techniques of \cite{JS3}.

\end{document}